\documentclass[letterpaper, 10 pt, conference]{ieeeconf}
\usepackage[utf8]{inputenc}

\IEEEoverridecommandlockouts    
\usepackage{cite}
\usepackage{amsmath,amssymb,amsfonts,amstext}
\usepackage{graphicx}
\usepackage{textcomp}

\usepackage[utf8]{inputenc}
\usepackage{babel}
\usepackage{balance}
\usepackage{comment}
\usepackage[position=bottom]{subfig}

\usepackage{todonotes}

\usepackage{mathtools}
\usepackage{float}

\usepackage{amsthm}

\newtheorem{remark}{Remark}

\newtheorem{problem}{Problem}

\usepackage{algorithm}
\usepackage{algorithmicx} 
\usepackage{algpseudocode} %

\makeatletter
\let\NAT@parse\undefined
\makeatother
\usepackage{hyperref}
\hypersetup{
  colorlinks=true,
    linkcolor= blue,
    allcolors=blue,
    citecolor = blue,
    filecolor=black,      
    urlcolor=blue,
    }

\title{\bf A Hierarchical Approach to Optimal Flow-Based Routing and Coordination of Connected and Automated Vehicles}

\author{Heeseung Bang$^{1,2}$, \textit{Student Member, IEEE}, Andreas A. Malikopoulos$^{2}$, \textit{Senior Member, IEEE}
    \thanks{This research was supported by NSF under Grants CNS-2149520 and CMMI-2219761.}
    \thanks{$^{1}$Department of Mechanical Engineering, University of Delaware, Newark, DE 19716, USA.} 
    \thanks{$^{2}$School of Civil and Environmental Engineering, Cornell University, Ithaca, NY 14850, USA. {\tt\small Email: \{hb489,amaliko\}@cornell.edu}}
}

\date{July 2022}

\begin{document}

\maketitle
\begin{abstract}
This paper addresses the challenge of generating optimal vehicle flow at the macroscopic level. Although several studies have focused on optimizing vehicle flow, little attention has been given to ensuring it can be practically achieved. To overcome this issue, we propose a route-recovery and eco-driving strategy for connected and automated vehicles (CAVs) that guarantees optimal flow generation. Our approach involves identifying the optimal vehicle flow that minimizes total travel time, given the constant travel demands in urban areas. We then develop a heuristic route-recovery algorithm to assign routes to CAVs. Finally, we present an efficient coordination framework to minimize the energy consumption of CAVs while safely crossing intersections. The proposed method can effectively generate optimal vehicle flow and potentially reduce travel time and energy consumption in urban areas.
\end{abstract}


\section{Introduction}

With increasing population and urbanization, traffic congestion has become a major issue in urban areas and transportation systems worldwide.
According to the Texas A\&M Transportation Institute \cite{Schrank2021}, in 2019, commuters in the US spent an extra $8.7$ billion hours and $3.5$ billion gallons of fuel only due to traffic congestion, which corresponds to \$$190$ billion of economic loss.

Connected and automated vehicles (CAVs) have attracted considerable attention since they offer the potential to reduce traffic congestion and fuel consumption by improving the efficiency of transportation systems \cite{Rios-Torres:2017aa}.
For example, numerous studies have addressed time and energy-efficient coordination problems in different traffic scenarios, i.e., merging roadways \cite{liu2023safety}, lane changes \cite{duan2023cooperative}, roundabouts \cite{alighanbari2020multi}, signal-free intersections \cite{Au2015,Malikopoulos2020}, and corridors \cite{lee2013,chalaki2020TITS,mahbub2020decentralized}.


Several efforts have been reported in the literature on routing single vehicles that optimize for different objectives \cite{de2017bi,houshmand2021combined,boewing2020vehicle,gammelli2021graph}.
These efforts focused more on selecting the best route for each vehicle regarding the events in the surrounding environment.
However, the resulting solutions may not be system optimal for multiple vehicles because they cannot consider the influence of other CAVs due to the demanding computation. 
Some studies followed sequential order \cite{Bang2022combined} or sought person-by-person optimal solutions \cite{Malikopoulos2021} among the CAVs \cite{Bang2022rerouting} for routing multiple CAVs while considering the influence of others. However, these approaches require further computational improvement to facilitate large-scale deployment.

Some other efforts have been made on routing CAVs, considering interactions between them with fewer computational issues.
Chu et al. \cite{chu2017dynamic} presented a routing method using a dynamic lane reversal approach. To consider the influence of other CAVs, they estimated travel time proportional to the number of CAVs on the roads.
Similarly, Mostafizi et al. \cite{mostafizi2021decentralized} designed a heuristic algorithm for a decentralized routing framework, which estimates the travel speeds inversely proportional to the number of CAVs on the same road.
As another way to avoid computational issues, a series of papers addressed the decentralized routing problem using the ant colony optimization technique, which uses a metric representing the traffic condition and updates the information in real-time \cite{bui2019aco,nguyen2021multiple,nguyen2021ant}.
However, in most of these methods, the estimations barely capture the actual traffic conditions. In addition, the methods are likely to yield locally optimal solutions due to the nature of the decentralized framework \cite{Malikopoulos2021}.

Many studies have also focused on addressing vehicle flow optimization problems at the macroscopic level to minimize travel time, energy consumption, and emissions.
Those studies considered flow optimization in different situations, e.g., considering relocating vehicles \cite{houshmand2019penetration,salazar2019congestion}, mixed traffic environment \cite{wollenstein2020congestion,wollenstein2021routing}, and charging constraints \cite{bang2021AEMoD}.
Most of these studies used a travel latency function to capture the influence of other vehicles.
Thus, they solved the optimization problem in a continuous domain, which enabled the use of the method for a large-scale network.
However, little attention has been given to ensuring that this optimal flow can be practically achieved from actual CAV movements.

To address all these challenges, we propose a hierarchical approach to optimal flow-based routing and coordination of CAVs. 
Our approach involves identifying the optimal vehicle flow that minimizes total travel time, given the constant travel demands in urban areas. We then develop a heuristic route-recovery algorithm to assign routes to CAVs that satisfy all travel demands while maintaining the optimal flow. At the lower level, we provide a coordination framework for CAVs to arrive at each road segment at their desired arrival time based on their assigned route and the desired flow. 
The main contributions of this work are as follows: (1) providing insights into combining different levels of routing methods, (2) developing a route-recovery algorithm that assigns routes to CAVs based on the flow, and (3) presenting a coordination framework that ensures the practical feasibility of optimal flow generation via CAV controls.

The remainder of this paper is organized as follows.
In Section \ref{sec:flow_optimization}, we present the traffic flow optimization in a macroscopic perspective and the route-recovery algorithm at a microscopic level.
In Section \ref{sec:coordination}, we provide the coordination framework at an intersection, and in Section \ref{sec:simulation}, we present simulation results to evaluate the performance of the proposed hierarchical approach.
Finally, in Section \ref{sec:conclusion}, we draw concluding remarks and discuss future research directions.

\section{Flow-Based Routing Framework}   \label{sec:flow_optimization}

In this section, we introduce a routing method based on optimal traffic flow.
In general, vehicle flow is captured by the number of vehicles passing through a point for a unit of time, and the traffic flow is analyzed and optimized from a macroscopic perspective without considering the movement of each vehicle.
Thus, we optimize the traffic flow at a macroscopic level and generate routes for each vehicle based on the optimal flow. Then, we provide a strategy for CAVs departing at depots to generate corresponding flow.

\subsection{Traffic Flow Optimization}

We consider an urban grid road network given by a directed graph $\mathcal{G} = (\mathcal{V},\mathcal{E})$, where $\mathcal{V}\subset\mathbb{N}$ is a set of vertices and $\mathcal{E}\subset\mathcal{V}\times\mathcal{V}$ is a set of roads.
Each vertex is either an intersection $r\in\mathcal{R}$ or a depot $q\in\mathcal{D}$, where $\mathcal{R}$ and $\mathcal{D}$ are the sets of intersections and depots, i.e., $\mathcal{V} = \mathcal{R} \cup \mathcal{D}$.
The depots are located at the entry/exit of intersections, implying that a CAV passing through a depot will exit one intersection and enter another.
Next, we determine travel information.
Let $M\in\mathbb{N}$ denote the total number of travel demands. For each travel demand $ m\in\mathcal{M} = \{1,\dots,M\}$, let $o_m, d_m\in\mathcal{D}$, and $\alpha_m\in\mathbb{R}_{>0}$ denote origin, destination, and demand rate, respectively. The demand rate is defined as the number of customers per unit of time.
For each travel demand $m\in\mathcal{M}$, $x_{ij}^m$ denotes the flow of CAVs traveling on the road $(i,j)$.
Then, the total flow on the road $(i,j)\in\mathcal{E}$ becomes $x_{ij} = \sum_{m\in\mathcal{M}} x_{ij}^m$.

To meet all the travel demands, the flow must satisfy the following constraints:
\begin{align}
    &\sum_{i:(i,j)\in\mathcal{E}}x_{ij}^m = \sum_{k:(j,k)\in\mathcal{E}}x_{jk}^m,~\forall m\in\mathcal{M},j\in\mathcal{V}\setminus\{o_m,d_m\},\label{eqn:con_x1}\\
    & \sum_{k:(j,k)\in\mathcal{E}}x_{jk}^m = \alpha_m,~\forall m\in\mathcal{M},j=o_m,\label{eqn:con_x2}\\
    & \sum_{i:(i,j)\in\mathcal{E}}x_{ij}^m = \alpha_m,~\forall m\in\mathcal{M},j=d_m.\label{eqn:con_x3}
\end{align}
Constraint \eqref{eqn:con_x1} guarantees that the incoming flow and the outgoing flow are the same at each node $j$, while constraints \eqref{eqn:con_x2} and \eqref{eqn:con_x3} match the flow at the origin and the destination with the demand rate $\alpha_m$. These constraints guarantee that the flow starts and ends at $o_m$ and $d_m$, respectively.

To estimate travel time on the roads, we use a travel time function proposed by the \textit{U.S. Bureau of Public Roads} (BPR) \cite{us1964traffic}.
By using the BPR function, the travel time on the road $(i,j)$ can be expressed as
\begin{equation}
    t_{ij}(x_{ij}) = t_{ij}^0\left(1+0.15\left(\frac{x_{ij}}{\gamma_{ij}}\right)^4\right), \label{eqn:BPR}
\end{equation}
where $t_{ij}^0 \in\mathbb{R}_{>0}$ is a free-flow travel time and $\gamma_{ij} \in\mathbb{R}_{>0}$ is capacity of the road $(i,j)\in\mathcal{E}$.
Next, we introduce a flow optimization problem using the BPR function.

\begin{problem}[Flow-based routing] \label{prb:flow_routing}
We find the optimal flow of CAVs by solving the following optimization problem:
\begin{equation}
\begin{aligned}
    \min_{\mathbf{x}} ~&J(\mathbf{x}) = \sum_{(i,j)\in\mathcal{E}} \bigg\{ t_{ij}(x_{ij})x_{ij}\bigg\} \label{eqn:flow_routing}\\
    \text{s.t. } & \eqref{eqn:con_x1} \text{--} \eqref{eqn:con_x3}, \mathbf{x} \geq 0,
\end{aligned}    
\end{equation}
where $\mathbf{x}$ is a vector of $x_{ij}$ for all $(i,j)\in\mathcal{E}$.
\end{problem}

Since BPR function \eqref{eqn:BPR} is convex in its domain and all the constraints are linear, Problem \ref{prb:flow_routing} becomes a convex optimization problem.
Problem \ref{prb:flow_routing} aims to minimize the total travel time for all vehicle flows.
Since Problem \ref{prb:flow_routing} is a strictly convex problem, it is guaranteed to have a unique global optimal solution.

\begin{remark}
In this paper, we focus on bridging the gap between flow optimization and coordination methods.
Hence, we consider $100\%$ penetration rate of CAVs in the network.
Nonetheless, it is possible to consider relocating trips \cite{salazar2019congestion}, mixed traffic for CAVs \cite{wollenstein2020congestion,wollenstein2021routing}, and even charging constraints \cite{bang2021AEMoD}.
In these cases, the optimization problem becomes a non-convex problem, which can be convexified using piece-wise affine functions to approximate the BPR function.
\end{remark}

\subsection{Route Recovery} \label{subsec:route_recovery}

The solution to Problem \ref{prb:flow_routing} determines the optimal flow at the roads for each travel demand $m$. Therefore, there could be multiple routes for CAVs to travel from origin to destination.
To assign a specific route to each CAV, we recover all different routes with corresponding flows well connected from origin to destination and match the demand rate.

Let $\mathcal{P}^m$ denote a set of routes for travel demand $m\in\mathcal{M}$, where $L^m\in\mathbb{N}$ denotes the number of routes, i.e., $L^m=|\mathcal{P}^m|$.
Each route $\mathcal{P}^m_l\in\mathcal{P}^m$, $l\in\{1,\dots,L^m\}$, is a tuple of the roads $(i,j)\in\mathcal{E}$ connected from $o_m$ to $d_m$, and $f_l^m$ is the corresponding flow.
We recover all routes based on the optimal flow using Algorithm \ref{Alg:route_recovery}.
The main idea of Algorithm \ref{Alg:route_recovery} is to search for connected flows from origin to destination.
During the search process, we prioritize checking the straight path to minimize the number of turns.
If we find the connected road, we update the flow with the smaller value among the current flow and the flow of the following road segment. This is because we are recovering routes, and the corresponding flow of the route must be the same for all road segments. We subtract the corresponding flow once it reaches the destination to avoid exploring the same route. We repeat the process until we find the routes for all flows.

\begin{algorithm}[t]
 \caption{Route Recovery}

\begin{algorithmic}[1] 
\For{$m\in\mathcal{M}$}
    \State{$l\gets 0$}
    \While{\texttt{not done}}
        \State{$l\gets l+1$;~~$i\gets o_m$;~~$f^m_l \gets$ $\alpha_m$}
        \While{$i \neq d_m$}
            \State{$j \gets$ next node} \Comment{Checking straight path first}
            \If{ $x_{ij} \neq 0 $}
                \State{$f^m_l \gets \min\{f^m_l,x_{ij}\}$}
                \State{$\mathcal{P}^m_l \gets$ Add road $(i,j)$}
                \State{$i \gets j$}
            \EndIf
        \EndWhile
        \For{$(i,j)\in\mathcal{P}^m_l$}
            \State{$x_{ij} \gets x_{ij} - f^m_l$}
        \EndFor
        \If{$\sum_{(i,j)}x_{ij} = 0$}
            \State{$L^m \gets l$}
            \State{\texttt{done} $\gets$ \texttt{True}}
        \EndIf
    \EndWhile
\EndFor
\end{algorithmic} \label{Alg:route_recovery}
\end{algorithm}

\begin{figure}
    \centering
    \includegraphics[width=0.6\linewidth]{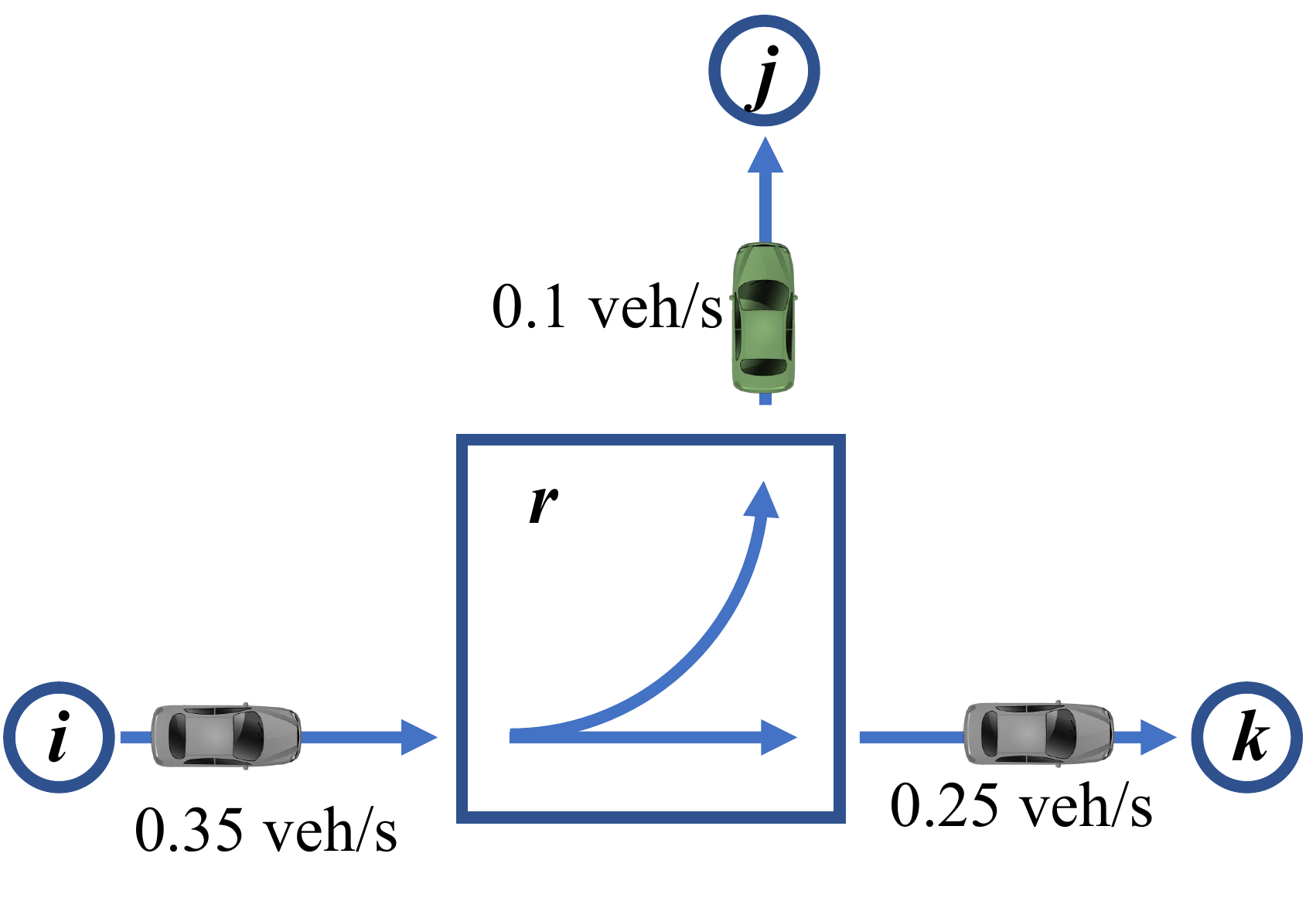}
    \caption{Illustration of CAVs departing at node $i\in\mathcal{D}$ for different routes. One is the roads $(i,r)$ and $(r,k)$, and the other is passing through $(i,r)$ and $(r,k)$.}
    \label{fig:departure_concept}
\end{figure}

\begin{figure}
    \centering
    \subfloat[Actual departure timing for each route]{\includegraphics[width=\linewidth]{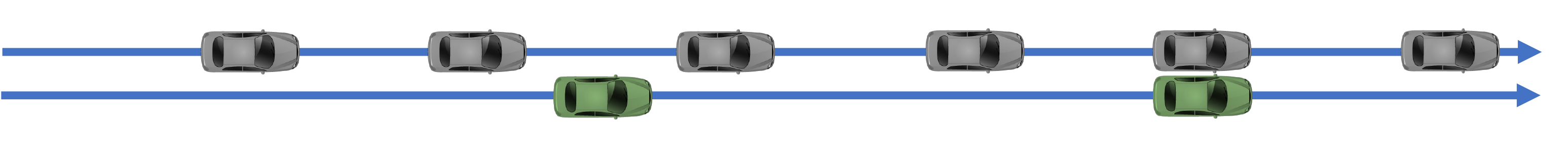}\label{b}}\\
    \subfloat[Possible departure timing at a depot]{\includegraphics[width=\linewidth]{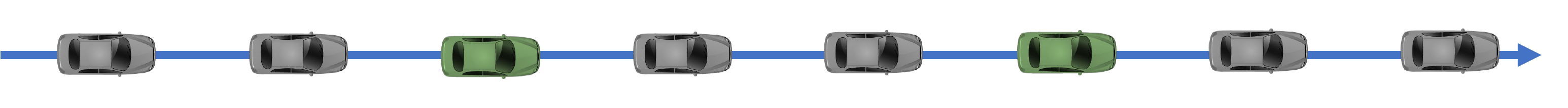}\label{c}}
    \caption{Illustration of departure time of CAVs for different flows.}
    \label{fig:departure_timing}
\end{figure}

Next, we discuss the rule for CAVs departing at the origin of the travel.
Suppose we obtain the optimal flow for two paths (Fig. \ref{fig:departure_concept}).
One group of CAVs generates the flow of $0.1$ veh/s while the other forms the flow of $0.25$ veh/s.
This implies that CAVs should enter/exit the roads $(i,r)$ and $(r,j)$ every $10$ seconds, and $(i,r)$ and $(r,k)$ every $4$ seconds.
The actual timing for CAVs to depart from depot $i$ to establish these flows can be illustrated in Fig. \ref{fig:departure_timing}-(a). This figure shows overlapping CAVs exist when they enter/exit at different frequencies. 
Therefore, we let CAVs depart at a uniform frequency depending on the optimal flow at the entrance road $(i,r)\in\mathcal{E}$ while maintaining the order of actual departure timing.
Next, we determine the actual time of each CAV $n$ crossing intersections.
Let $t_n^0$ be the entry time of CAV $n$ at an intersection.
Then, we select the exit time based on the following conditions:
\begin{equation}
    t_n^f = \max\left\{t_n^0+t_{ij}(x^*_{ij})+t_{jk}(x^*_{jk}),~t_p^f+\frac{1}{x^*_{jk}}\right\}. \label{eqn:terminal_time}
\end{equation}
Here, $t_p^f$ is the exit time of preceding CAV $p$ that exits right before CAV $n$, and $x^*_{(\cdot)}$ is the optimal flow on the given road.
We select the maximum value of two options, where the first possible value in \eqref{eqn:terminal_time} is the estimated arrival time at the exit node using the BPR function.
If no preceding CAV exits to the same node, CAV $n$ moves along the estimated travel time.
However, the entering frequencies for different entry nodes are not synchronized.
Thus, a CAV $p$ with overlapping exit time may exist from different entry nodes.
To avoid this situation, we let CAV $n$ select the second possible value in \eqref{eqn:terminal_time} so that CAV $n$ can exit after CAV $p$ while maintaining the optimal frequency.

\section{Flow-Matching Coordination at Intersections} \label{sec:coordination}

In this section, we present a coordination framework for CAVs to cross a signal-free intersection.
Given the routes for all CAVs and specific entry/exit times for all the intersections, CAVs plan their trajectory to cross an intersection while reducing energy consumption and avoiding collisions.

\begin{figure}
    \centering
    \includegraphics[width=0.65\linewidth]{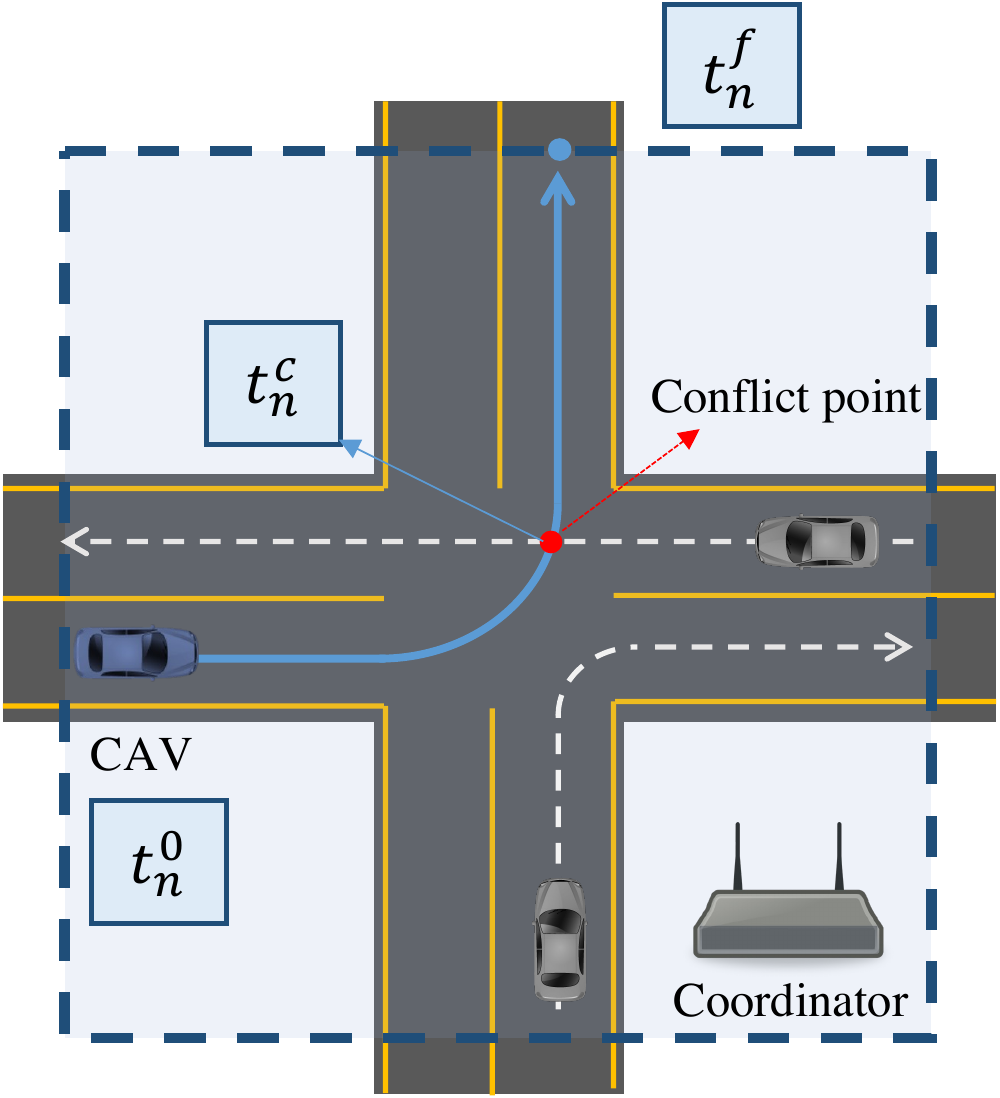}
    \caption{Coordination of CAVs at a signal-free intersection.}
    \label{fig:intersection}
\end{figure}

We consider a single-lane intersection as illustrated in Fig. \ref{fig:intersection}. It consists of four entries and four exits, and CAVs can take twelve different paths in total, as we neglect U-turns.
At each intersection $r\in\mathcal{R}$, we have a \textit{coordinator} that communicates with CAVs and provides geometry information and trajectories of  other CAVs.
When CAVs enter the intersection, they plan their trajectories with respect to the existing CAVs in the intersection. This means that CAVs fix their trajectories in the order of entrance.
There is a set of points $\mathcal{C}\subset\mathbb{N}$ called \textit{conflict points} at which the lateral collision may occur among the CAVs.

Let $\mathcal{N}_r(t)$ be the set of CAVs in the intersection $r\in\mathcal{R}$ at time $t\in\mathbb{R}_{\geq 0}$.
For each CAV $n\in\mathcal{N}_r(t)$, we consider double integrator dynamics
\begin{equation}
    \begin{aligned}
    \dot{s}_n(t) &= v_n(t),\\
    \dot{v}_n(t) &= u_n(t), \label{eqn:dynamics}
    \end{aligned}
\end{equation}
where $s(t)\in\mathbb{R}_{\geq0}$ is the distance between the entry point to the current location, $v(t)\in\mathbb{R}_{\geq0}$ denotes a speed, and $u(t)\in\mathbb{R}$ denotes control input of CAV $n$.
For each CAV $n$, we impose state/input constraints as
\begin{align}
    u_{n,\text{min}} & \leq u_n(t) \leq u_{n,\text{max}}, \label{eqn:ulim}\\
    0 < v_{\text{min}} & \leq v_n(t) \leq v_{\text{max}}, \label{eqn:vlim}
\end{align}
where $u_{n,\text{min}}$, $u_{n,\text{max}}$ are the minimum and maximum control inputs and $v_{\text{min}}$, $v_{\text{max}}$ are the minimum and maximum allowable speeds, respectively.

Next, we introduce safety constraints.
To ensure rear-end safety, we impose constraints on the minimum distance between preceding CAV $p\in\mathcal{N}_r(t)$ and following CAV $n\in\mathcal{N}_r(t)$, i.e.,
\begin{equation}
    s_p(t) - s_n(t) \geq \delta, \label{eqn:rear-end}
\end{equation}
where $\delta \in\mathbb{R}_{>0}$ is a safety distance.
For lateral safety, suppose there exists CAV $p$ that shares a conflict point with CAV $n$.
Then, we impose the following constraint,
\begin{equation}
    |t_p^c - t_n^c | \geq \tau, \label{eqn:lateral}
\end{equation}
where $\tau\in\mathbb{R}_{\geq0}$ is a safety time headway, and $t_p^c$, $t_n^c$ are the arrival time of CAV $p$ and CAV $n$ at a conflict point $c\in\mathcal{C}$, respectively.

In our coordination framework, we minimize the speed transition of CAVs to conserve momentum and reduce fuel consumption.
\begin{problem}
To find minimum-energy trajectory, each CAV $n\in\mathcal{N}_r(t)$ solves the following optimization problem
\begin{align} \label{eqn:energy-optimal}
    \min &\text{~} \frac{1}{2} \int_{t_n^0}^{t_n^f} u_n^2(t) dt\\
    \emph{s.t. }& \eqref{eqn:dynamics} -  \eqref{eqn:lateral} \notag.
\end{align}
where $t_n^0$ is an entry time, $t_n^f$ is an exit time of CAV $n$, respectively.
    \label{prb:constrained_energy_optimal}
\end{problem}

To address Problem \ref{prb:constrained_energy_optimal}, we use different approaches depending on the activation of safety constraints.

\textbf{\textit{1) No safety constraints violated}}: In case of no safety constraints violation, we use a solution to the unconstrained optimization problem. 
The unconstrained energy-optimal trajectory is \cite{chalaki2021CSM}
\begin{align} \label{eq:optimalTrajectory}
    u_n(t) &= 6 a_n t + 2 b_n, \notag \\
    v_n(t) &= 3 a_n t^2 + 2 b_n t + c_n, \\
    s_n(t) &= a_n t^3 + b_n t^2 + c_n t + d_n, \notag
\end{align}
where $a_n, b_n, c_n$, and $d_n$ are constants of integration. The constrained optimal solution of \eqref{eq:optimalTrajectory} is reported in \cite{malikopoulos2019ACC}.
We obtain these constants from the following boundary conditions
\begin{align}
     s_n(t_n^0) &= 0,\quad  v_n(t_n^0)= v_n^0 , \label{eq:bci}\\
     s_n(t_n^f)&=s_n^f,\quad v_n(t_n^f)=v_n^f. \label{eq:bcf}
\end{align}
Here, $s_n^f$ is the road length from the entry to the exit, $v_n^0$ is the entry speed, and $v_n^f$ is the exit speed of CAV $n$.

Next, we verify if \eqref{eq:optimalTrajectory} satisfies speed limits and control input constraints for $t\in[t_n^0,t_n^f]$. Note that $v_n^f$ has not been determined yet. Therefore, we can select $v_n^f$ that allows the trajectory to satisfy state/input constraints by solving the following optimization problem. 

\begin{problem} \label{prb:exit_speed}
Each CAV $n\in\mathcal{N}_r(t)$ selects its exit speed by solving the following optimization problem:
\begin{equation}
\begin{aligned}
    \min & \left(v_n^f-\bar{v}\right)^2\\ 
    \text{s.t. } & \eqref{eqn:ulim}, \eqref{eqn:vlim}, \eqref{eq:optimalTrajectory} - \eqref{eq:bcf},
\end{aligned}    
\end{equation}
where $\bar{v}$ is an estimated average speed in the next road segment.
\end{problem}

The estimated average speed $\bar{v}\in\mathbb{R}_{>0}$ can be simply selected using the estimated travel time, i.e., $\bar{v}=S/t_{ij}(x^*_{ij})$, where $S\in\mathbb{R}_{>0}$ and $x^*_{ij}$ are the length and the optimal flow on the next road segment $(i,j)\in\mathcal{E}$, respectively. If the trajectory does not violate the state/input constraints, the exit speed will be $v_n^f=\bar{v}$.

\textbf{\textit{2) Lateral safety violated}}: If the unconstrained energy-optimal trajectory violates the lateral safety constraints, we could use the method described in \cite{Malikopoulos2020} to solve Problem \ref{prb:constrained_energy_optimal}.
This method suggests that we piece the constrained and unconstrained arcs together, but it naturally yields difficulties in implementing real-time solvers.
To avoid this issue, we search for a way-point that guarantees the trajectories not to violate lateral safety constraints. Let $\mathcal{W}_n$ be a set of possible way-points for CAV $n\in\mathcal{N}_r(t)$.
We select possible waypoints as the boundaries of \eqref{eqn:lateral} to minimize the time gap and maximize throughput.
For example, one possible way-point for CAV $n$ is $(t_p^c+\tau,s_n^c)$, which implies that CAV $n$ passes the location of the conflict point $s_n^c$ with a safety time gap $\tau$ after CAV $p$ passes.

We piece two unconstrained energy-optimal trajectories, one from the entry to the way-point and another from the way-point to the exit. Note that both trajectories' state/input constraints need to be considered in selecting speed at the way-point $v_n^w, w\in\mathcal{W}$.
By setting the speed $v_n^w$ as exit/entry speeds for those two trajectories, we ensure that CAV $n$ has a continuous speed profile and can compute its trajectory using the same procedure as in the unconstrained case.
This approach sacrifices some of the optimality, but each piece of trajectory is still energy-optimal for given boundary conditions, and it provides a feasible solution with simple computations.

\textbf{\textit{3) Rear-end safety violated}}: CAVs depart from the origin with uniform time gaps and also exit the intersection with guaranteed time gap from \eqref{eqn:terminal_time}.
Thus, rear-end safety is first violated when the preceding CAV $p$ slows down to avoid a lateral collision.
Then, the following CAV $n$ requires a new way-point that does not violate rear-end safety.
We let CAV $n$ follow the way-point of its preceding CAV $p$ while keeping the safety distance.
For example, CAV $n$ can select $(t_p^c,s_n^c-\delta)$ as a way-point, where $c\in\mathcal{C}$ is a conflict point at which CAV $p$ needs to avoid a lateral collision.
We set the speed of CAV $n$ at this point to be the same as CAV $p$ so that the rear-end safety constraint is not violated near the point.

\begin{figure*}
    \centering
    \subfloat[Total travel demands]{\includegraphics[width=0.33\linewidth]{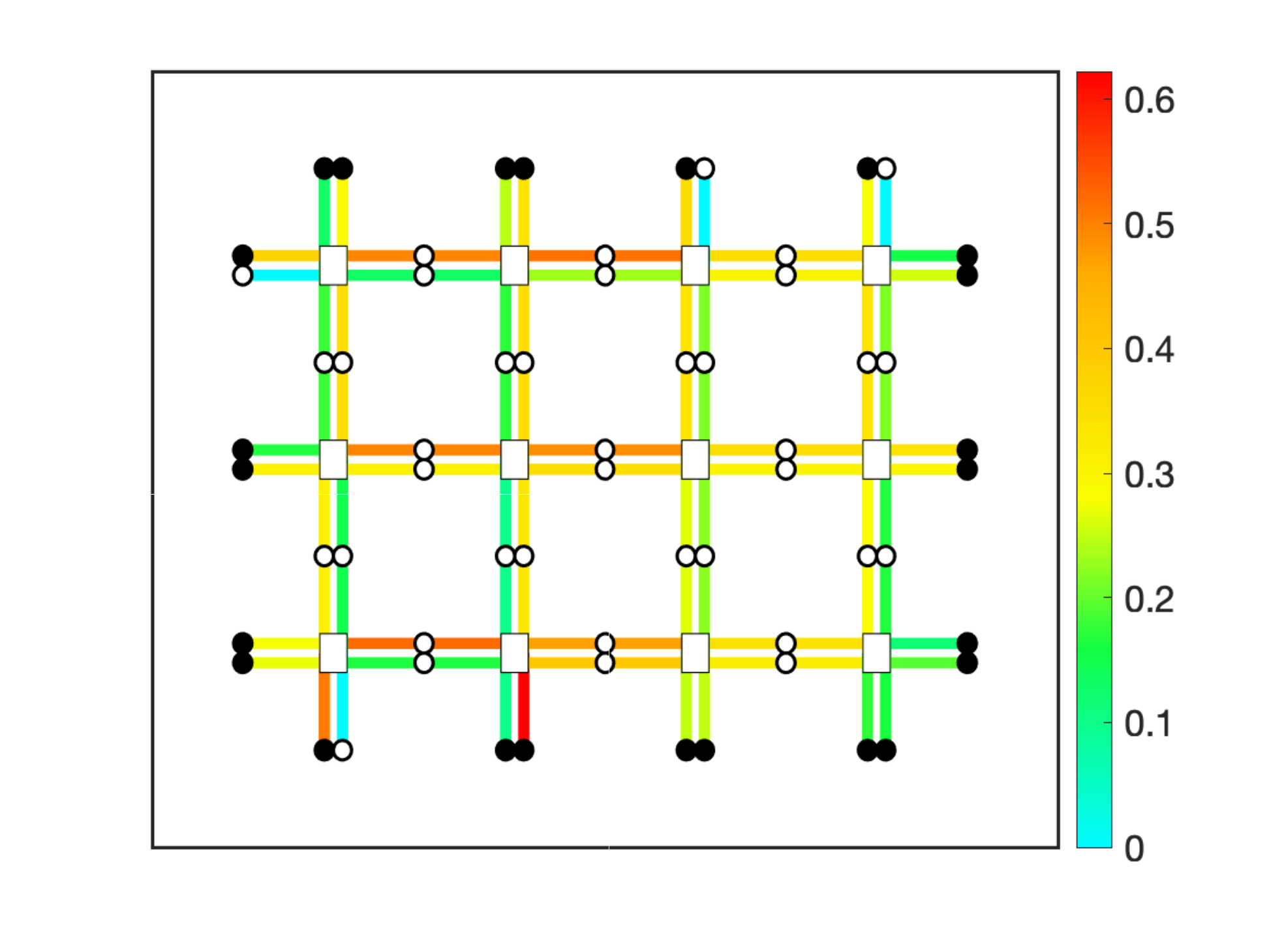}}
    \subfloat[Selected travel demand \#1]{\includegraphics[width=0.33\linewidth]{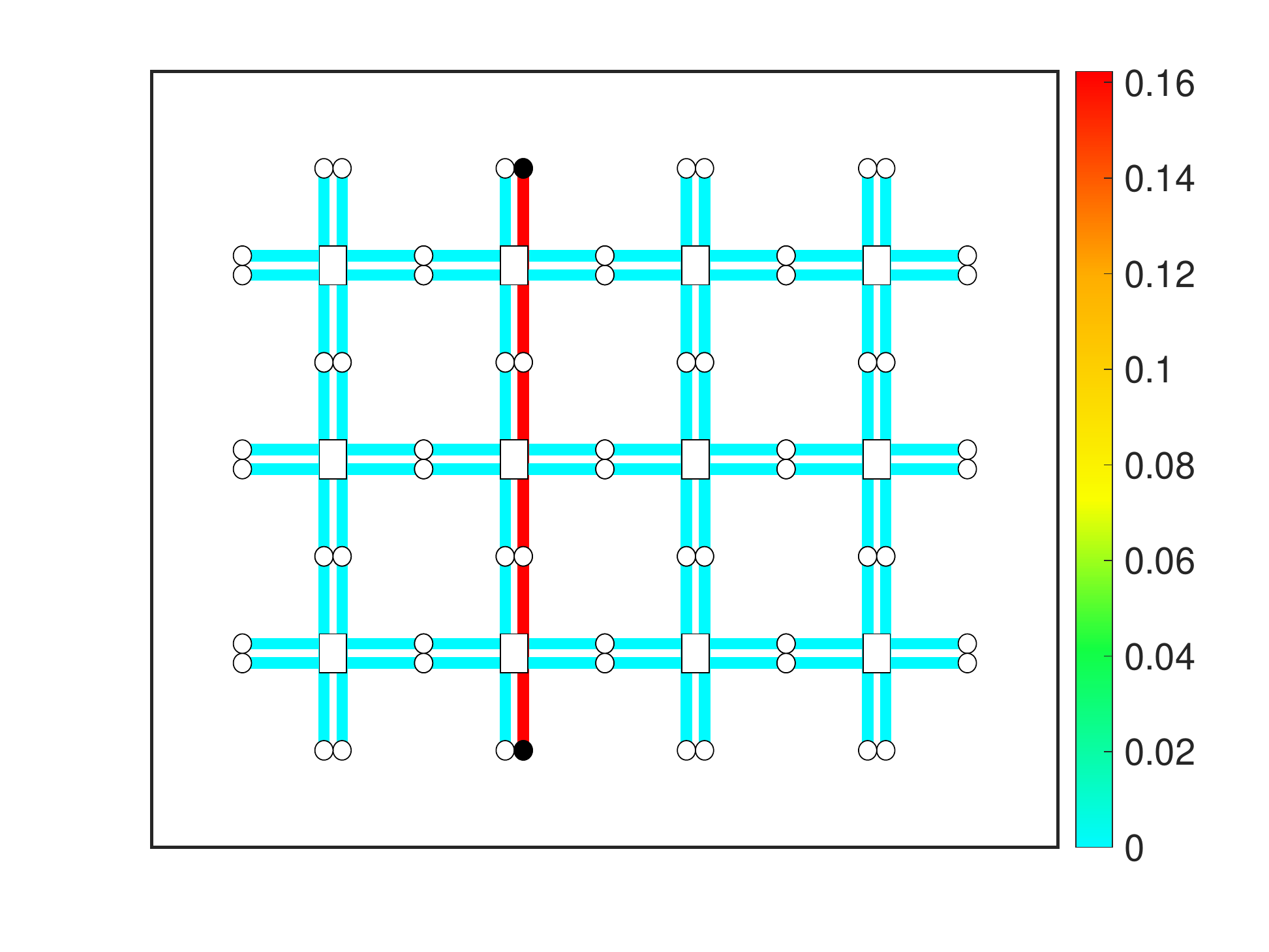}}
    \subfloat[Selected travel demand \#2]{\includegraphics[width=0.33\linewidth]{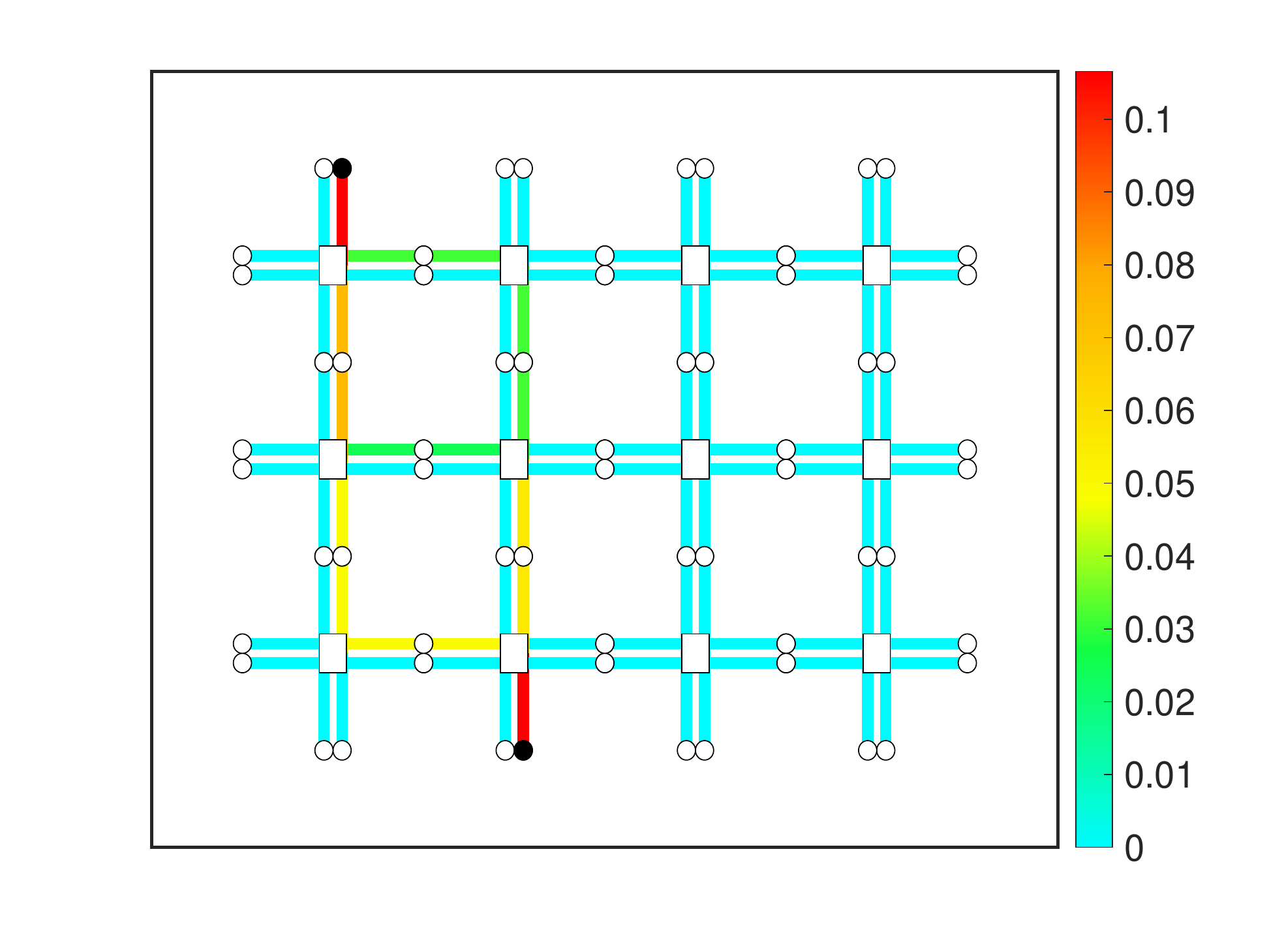}}
    \caption{Graphical illustration of optimal flow on a grid road network.}
    \label{fig:total_flow}
\end{figure*}

\section{Simulation Results} \label{sec:simulation}

This section provides the simulation results at different levels to analyze whether our method is implementable in macroscopic and microscopic perspectives. 
First, for the routing level, we present the optimal flow for different travel demands at a road network.
Then, we focus on a single intersection for coordination and show that CAVs can safely cross the intersection within a given travel time.

\subsection{Vehicle Flow Optimization}
We considered a road network with $12$ intersections, $62$ depots, and $96$ road segments.
The length of each road segment is $200$ m. 
We randomly generated $30$ travel demands with origin, destination, and demand rate information.
Figure \ref{fig:total_flow} illustrates the simulation results for (a) total travel demands and some (b), (c) selected travel demands. In the figure, all lines represent road segments where the color characterizes the flow, and the white squares exhibit intersections. All  circles are the depots, whereas black circles are the origin or destination of travel demands.
Although some congested roads exist due to uneven travel demands, results show that the flows are well distributed to minimize the total flow.
For example, Fig. \ref{fig:total_flow}-(b) is the optimal flow for a specific travel demand.
This is the optimal flow and route for such an origin-destination pair because if CAVs in this flow turn left or right, they increase their own travel time and other CAVs' travel time. After all, all the other roads are congested, as shown in Fig. \ref{fig:total_flow}-(a).
On the other hand, another flow illustrated in Fig. \ref{fig:total_flow}-(c) was separated into $3$ different routes not to increase flow on the congested roads as much as possible.


\subsection{Coordination at an intersection}
To coordinate CAVs, we selected a specific intersection (the bottom left one in the network).
We considered $100$ CAVs and generated entry/exit time based on the method described in Section \ref{subsec:route_recovery}.


\begin{figure*}
    \centering
    \includegraphics[width=0.95\linewidth]{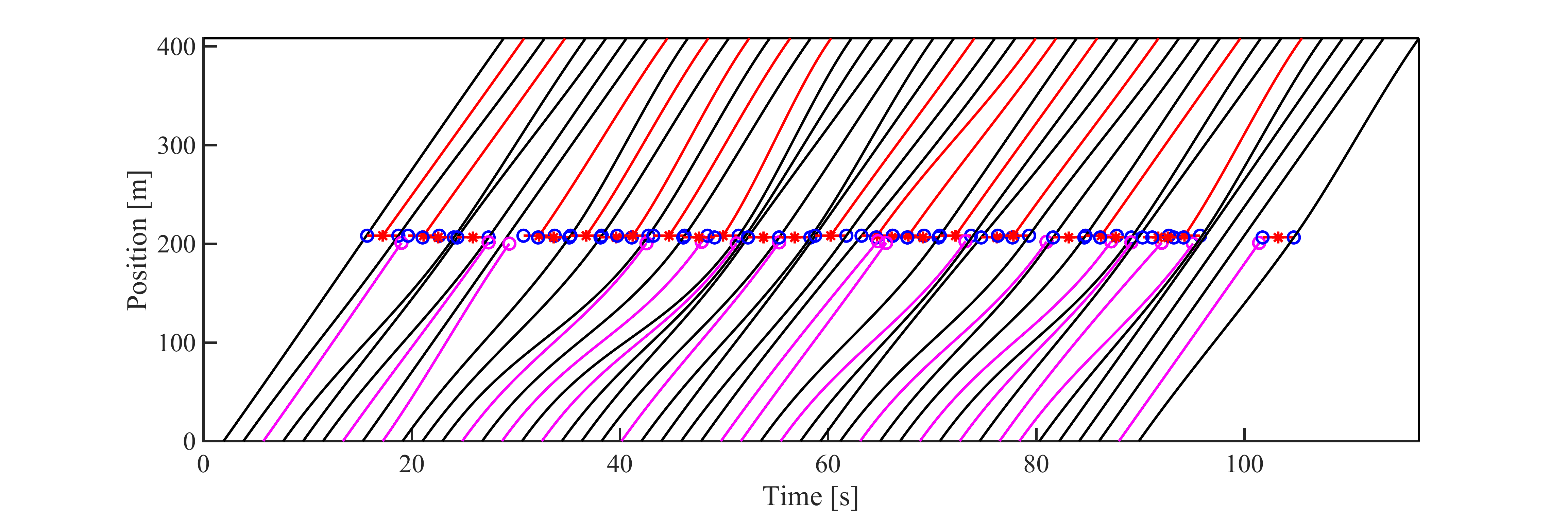}
    \caption{Position trajectories of CAVs using specific paths.}
    \label{fig:coordination}
\end{figure*}

We conducted simulations for $100$ CAVs, and the results are presented in Fig. \ref{fig:coordination}.
The graph shows all the trajectories of CAVs passing through specific paths in an intersection.
The red and purple trajectories are the CAVs merging to and splitting from the given path.
As shown in the figure, we efficiently assigned CAVs to cross the intersection by letting them pass through the boundaries of safety constraints.
The results also explain that all the CAVs selected proper exit time based on availability and that there are no cascading delay effects for both widely/densely distributed CAVs.
In the case of piecing two unconstrained trajectories, we set $\bar{v}$ to be the average of entry and exit speeds.
In the simulation, Solution to Problem \ref{prb:exit_speed} was always $\bar{v}$, which implies that state/input constraints were never violated with the given boundary conditions.


\section{Concluding Remarks} \label{sec:conclusion}

In this paper, we proposed a hierarchical approach for generating optimal vehicle flow in urban areas using a route-recovery algorithm and a coordination framework for CAVs.
The proposed approach minimizes total travel time and ensures the practical feasibility of achieving optimal flow. Our coordination framework minimizes energy consumption and prevents collisions while crossing intersections.

Future work involves analyzing the results in the entire network and assessing the trade-off between optimality and computational cost when piecing two trajectories together in the coordination framework.
Another interesting direction for future research is to  investigate different travel time functions best suited for an urban road network.


\bibliographystyle{IEEEtran}
\bibliography{Bang, IDS}

\end{document}